\documentclass[11pt]{article}
\usepackage{amsthm}
\usepackage{pxfonts}
\usepackage{yfonts}
\usepackage{dsfont}
\usepackage{graphicx}
\usepackage{relsize}
\usepackage{color}

\def\bel{\begin{equation}\label}
\def\eeq{\end{equation}}
\def\ds{\displaystyle}

\def\mt{\longrightarrow}
\def\v{\vskip 1em}

\def\Rec{{\bf R}}
\def\R{\mathbb R}

\def\C{\mathfrak{B}}

\def\N{{\bf N}}

\def\S{{\bf S}}

\def\F{{\bf F}}

\def\P{{\bf P}}
\def\Q{{\bf Q}}

\def\L{{\bf L}}

\def\V{{\bf V}}

\def\T{{\bf T}}

\def\p{{\partial}}
\def\a{{\bf a}}

\def\bar{\overline}

\def\vol{{\bf vol}}

\def\I{{\bf I}}

\def\M{{\bf M}}

\def\alpha{\alphaup}
\def\beta{\betaup}
\def\gamma{\gammaup}
\def\delta{\deltaup}
\def\theta{\thetaup}
\def\xi{{\xiup}}
\def\eta{{\etaup}}
\def\tau{{\tauup}}
\def\rho{{\rhoup}}
\def\phi{{\phiup}}
\def\psi{{\psiup}}
\def\lambda{{\lambdaup}}
\def\omega{\omegaup}
\def\varphi{{\varphiup}}
\def\gamma{{\gammaup}}

\newtheorem{remark}{Remark}[section]

\begin{document}

\[\hbox{\LARGE{\bf Fractional integration associated with Zygmund dilations}}\]

\[\hbox{Chuhan Sun}\qquad\hbox{and}\qquad\hbox{Zipeng Wang}\]
\begin{abstract}
First, we study a family of fractional integral operator defined as
\[\I_{\alpha\beta\vartheta}f(u,v,t)~=~\iiint_{\R^{2n+1}} f(\xi,\eta,\tau)\V^{\alpha\beta\vartheta}[(u,v,t)\odot(\xi,\eta,\tau)^{-1}]d\xi d\eta d\tau\]
where $\odot$ denotes the multiplication law on a real Heisenberg group.
\v
$\V^{\alpha\beta\vartheta}$ is a distribution in $\R^{2n+1}$ satisfying Zygmund dilations.
A characterization is established between $\I_{\alpha\beta\vartheta}\colon\L^p(\R^{2n+1})\mt\L^q(\R^{2n+1})$ and necessary constraints consisting of $\alpha,\beta\in\R$ and $\vartheta\ge0$ for $1<p< q<\infty$. 
\v
Moreover, denote $\hbox{\bf R}=\Q\times\P\times I\subset\R^n\times\R^n\times\R$: $I$ is an interval and $\Q,\P$ are cubes parallel to the coordinates.  For $\alpha,\beta\in\R$, we consider 
\[\begin{array}{lr}\ds
\M_{\alpha\beta}f(u,v,t)~=~
\\ \ds
\sup_{\Rec\ni(0,0,0) \colon \vol \{I\}=\vol\{\Q\}^{1\over n}\vol\{\P\}^{1\over n}} \vol \{\Q\}^{{\alpha\over n}-1} \vol \{\P\}^{{\alpha\over n}-1} \vol \{I\}^{\beta-1}

\iiint_\Rec \left|f [(u,v,t)\odot(\xi,\eta,\tau)^{-1}]\right|d\xi d\eta d\tau.
\end{array}
\]
We show $\M_{\alpha\beta}\colon\L^p(\R^{2n+1})\mt\L^q(\R^{2n+1})$ for $1<p\leq q<\infty$ whenever ${\alpha+\beta\over n+1}={1\over p}-{1\over q}$.
\end{abstract}
\section{Introduction}
\setcounter{equation}{0}
A fractional integral operator $\T_\a$ is initially defined on $\R^\N$ as
\bel{I_a}
\T_\a f(x)~=~\int_{\R^\N}f(y) \left[{1\over|x-y|}\right]^{\N-\a} dy,\qquad 0<\a<\N.
\eeq
In 1928, Hardy and Littlewood \cite{Hardy-Littlewood} have obtained an regularity theorem for $\T_\a$ when $\N=1$. Ten years later, Sobolev \cite{Sobolev} made extensions on every higher dimensional space.

$\diamond$ {\small Throughout, $\C>0$ is regarded as a generic constant depending on its sub-indices.}
\v
{\bf Hardy-Littlewood-Sobolev theorem}~~~{\it Let $\T_\a$ defined in (\ref{I_a}) for $0<\a<\N$. We have
\bel{HLS Ineq}
\begin{array}{cc}\ds
\left\| \T_\a f\right\|_{\L^q(\R^\N)}~\leq~\C_{p~q}~\left\| f\right\|_{\L^p(\R^\N)},\qquad 1<p<q<\infty
\\\\ \ds
\hbox{\small{if and only if}}\qquad {\a\over \N}~=~{1\over p}-{1\over q}.~~~~~~~
\end{array}
\eeq
}\\
This classical result was first re-investigated by Folland and Stein \cite{Folland-Stein} on Heisenberg groups.

In this paper, $\odot$ denotes a multiplication law:
\bel{multiplication law}
\begin{array}{cc}\ds
(u,v,t)\odot(\xi,\eta,\tau)~=~\Big[u+\xi, v+\eta,t+\tau+\mu\big(u\cdot\eta-v\cdot\xi\big)\Big],\qquad\mu\in\R
\end{array}
\eeq
for every $(u,v,t)\in\R^n\times\R^n\times\R$ and $(\xi,\eta,\tau)^{-1}=(-\xi,-\eta,-\tau)\in\R^n\times\R^n\times\R$.

Let $0<\delta<n+1$. Consider
\bel{T_a}
\S_\delta f(u,v,t)~=~\iiint_{\R^{2n+1}}f(\xi,\eta,\tau)\Omega^\delta\Big[(u,v,t)\odot(\xi,\eta,\tau)^{-1}\Big]d\xi d\eta d\tau.
\eeq
$\Omega^\delta$ is a distribution in $\R^{2n+1}$ agree with
\bel{Omega^a}
\Omega^\delta(u,v,t)~=~\left[{1\over |u|^2+|v|^2+|t|}\right]^{n+1-\delta},\qquad \hbox{\small{$(u,v,t)\neq(0,0,0)$}}.
\eeq

{\bf Folland-Stein theorem} ~~{\it Let $\S_\delta$ defined in (\ref{T_a})-(\ref{Omega^a}) for $0<\delta<n+1$. We have
\bel{Folland-Stein theorem}
\begin{array}{cc}\ds
\left\| \S_\delta f\right\|_{\L^q(\R^{n+1})}~\leq~\C_{p~q}~\left\| f\right\|_{\L^p(\R^{2n+1})},\qquad 1<p<q<\infty
\\\\ \ds
\hbox{if and only if}\qquad {\delta\over n+1}~=~{1\over p}-{1\over q}.
\end{array}
\eeq}\\
The best constant for the $\L^p\mt\L^q$-norm inequality in (\ref{Folland-Stein theorem}) is found by Frank and Lieb \cite{Frank-Lieb}. A discrete analogue of this result has been obtained by Pierce \cite{Pierce}. Recently, the regarding commutator estimates are established  by Fanelli and~Roncal \cite{Luca-Luz}.

Next, we give a multi-parameter extension to {\bf Folland-Stein theorem} by replacing $\Omega^\delta$ with a larger kernel having singularity on every coordinate subspace. 
From (\ref{Omega^a}), we have
\[
\Omega^\delta(u,v,t)~\leq~\Bigg[{1\over |u||v|+|t|}\Bigg]^{n+1-\delta},\qquad \hbox{\small{$(u,t)\neq(0,0)$~~ or~~ $(v,t)\neq(0,0)$}}.
\]
A direct computation shows
\[\begin{array}{lr}\ds
\left[{1\over |u||v|+|t|}\right]^{n+1-\delta}
~\approx~ \left[{1\over |u|^2|v|^2+t^2}\right]^{{n+1\over 2}-{\delta\over2}}
\\\\ \ds~~~~~~~~~~~~~~~~~~~~~~~~~~~~
~=~|u|^{{\delta\over 2}-{n+1\over 2}}|v|^{{\delta\over 2}-{n+1\over 2}}|t|^{{\delta\over 2}-{n+1\over 2}} \left[{|u||v||t|\over |u|^2|v|^2+t^2}\right]^{{n+1\over 2}-{\delta\over2}}
\\\\ \ds~~~~~~~~~~~~~~~~~~~~~~~~~~~~
~=~|u|^{\big[{\delta\over 2}+{n-1\over 2}\big]-n}|v|^{\big[{\delta\over 2}+{n-1\over 2}\big]-n}|t|^{\big[{\delta\over 2}-{n-1\over 2}\big]-1} \left[{|u||v|\over |t|}+{|t|\over|u||v|}\right]^{-\big[{n+1\over 2}-{\delta\over2}\big]}
\end{array}
\]
whenever  $(u,t)\neq(0,0)$ or $(v,t)\neq(0,0)$. 

Above estimates  lead us to the following assertion.
Let $\alpha, \beta\in\R$ and $\vartheta\ge0$. $\V^{\alpha\beta\vartheta}$ is a distribution in $\R^{2n+1}$ agree with 
\bel{V}\begin{array}{lr}\ds
\V^{\alpha\beta\vartheta}(u,v,t)~=~|u|^{\alpha-n}|v|^{\alpha-n}|t|^{\beta-1} \Bigg[ {|u||v|\over |t|}+{|t|\over |u||v|}\Bigg]^{-\vartheta},
\qquad\hbox{\small{$u\neq0, v\neq0$, $t\neq0$}}.
\end{array}
\eeq
\begin{remark}
By taking into account $\alpha={\delta\over 2}+{n-1\over 2}$, $\beta={\delta\over 2}-{n-1\over 2}$  and $\vartheta={n+1\over 2}-{\delta\over2}$
for $0<\delta<n+1$, we find $\alpha>n\beta$ and $\vartheta={n+1\over 2}-{\delta\over2}>{\alpha-n\beta\over n+1}$. 
\end{remark}
Define 
\bel{I}
\begin{array}{lr}\ds
\I_{\alpha\beta\vartheta} f(u,v,t)~=~\iiint_{\R^{2n+1}} f\left(\xi,\eta,\tau\right) \V^{\alpha\beta\vartheta}\Big[(u,v,t)\odot(\xi,\eta,\tau)^{-1}\Big]d\xi d\eta d\tau.
\end{array}
\eeq
Observe that
\bel{Zygmund dilation}
\V^{\alpha\beta\vartheta}\Big[(ru,sv,rst)\odot(r\xi,s\eta,rs\tau)^{-1}\Big]~=~r^{\alpha+\beta-n-1}s^{\alpha+\beta-n-1} \V^{\alpha\beta\vartheta}\Big[(u,v,t)\odot(\xi,\eta,\tau)^{-1}\Big]
\eeq
for every $ r,s>0$.

A convolution operator of this type is said to be associated with Zygmund dilations. 
Certain singular integral operators  defined on a Heisenberg group associated with Zygmund dilations have been systematically studied, for example  by Phong and Stein \cite{Phong-Stein}, Ricci and Stein \cite{Ricci-Stein} and M\"{u}ller, Ricci and Stein \cite{Muller-Ricci-Stein}. Much less is known  for fractional integrals in this direction.
\v
{\bf Theorem One}~~{\it Let $\I_{\alpha\beta\vartheta}$ defined in  (\ref{V})-(\ref{I}) for $\alpha,\beta\in\R$ and $\vartheta\ge0$. We have
\bel{Result One}
\begin{array}{cc}\ds
\left\| \I_{\alpha\beta\vartheta} f\right\|_{\L^q(\R^{2n+1})}~\leq~\C_{p~q}~\left\| f\right\|_{\L^p(\R^{2n+1})},\qquad 1<p<q<\infty
\\\\ \ds
\hbox{if and only if}\qquad \vartheta\ge{|\alpha-n\beta |\over n+1},\qquad
{\alpha+\beta\over n+1}~=~{1\over p}-{1\over q}.
\end{array}
\eeq}

Denote $\hbox{\bf R}=\Q\times\P\times I\subset\R^n\times\R^n\times\R$: $I$ is an interval and $\Q,\P$ are cubes parallel to the coordinates.  For $\alpha,\beta\in\R$, we define
\bel{M Z}
\begin{array}{lr}\ds
\M_{\alpha\beta}f(u,v,t)~=~
\\ \ds
\sup_{\Rec\ni(0,0,0) \colon \vol \{I\}=\vol\{\Q\}^{1\over n}\vol\{\P\}^{1\over n}} \vol \{\Q\}^{{\alpha\over n}-1} \vol \{\P\}^{{\alpha\over n}-1} \vol \{I\}^{\beta-1}
\iiint_\Rec \left|f [(u,v,t)\odot(\xi,\eta,\tau)^{-1}]\right|d\xi d\eta d\tau.
\end{array}
\eeq
For $\M_{\alpha\beta}$ defined in the Euclidean space, in particular at $\alpha=\beta=0$, the regarding $\L^p$- theorem and its weighted analogue have been well established. See the paper by Ricci and Stein \cite{Ricci-Stein} and Fefferman and Pipher \cite{Fefferman-Pipher}.
\v
{\bf Theorem Two}~~{\it Let $\M_{\alpha\beta}$ defined in  (\ref{M Z}) for $\alpha,\beta\in\R$. We have
\bel{Result Two}
\begin{array}{cc}\ds
\left\| \M_{\alpha\beta} f\right\|_{\L^q(\R^{2n+1})}~\leq~\C_{p~q}~\left\| f\right\|_{\L^p(\R^{2n+1})},\qquad 1<p\leq q<\infty
\\\\ \ds
\hbox{whenever}\qquad
{\alpha+\beta\over n+1}~=~{1\over p}-{1\over q}.
\end{array}
\eeq}
\begin{remark}  For $1<p=q<\infty$ and $\alpha+\beta=0$, the result in (\ref{Result Two}) is an implication from the $\L^p$-theorem of strong maximal functions defined in Heisenberg groups, established by M. Chirst \cite{Michael Christ 2}.
\end{remark}
Later, we shall find
\bel{M<M}
\M_{\alpha\beta}f(u,v,t)~\leq~\M_\gamma f(u,v,t),\qquad \gamma={\alpha+\beta\over n+1},\qquad 0\leq\gamma<1.
\eeq
$\M_\gamma, 0\leq\gamma<1$ is defined as
\bel{M_gamma}
\begin{array}{lr}\ds
\M_\gamma f(u,v,t)~=~\sup_{\Rec\ni(0,0,0)}
\vol \{\Rec\}^{\gamma-1} \iiint_\Rec\left|f [(u,v,t)\odot(\xi,\eta,\tau)^{-1}]\right|d\xi d\eta d\tau.
\end{array}
\eeq
\begin{remark}
The supremum in (\ref{M_gamma}) is taking over any rectangle  $\Rec\subset\R^{2n+1}$ parallel to the coordinates. 
\end{remark}
$\M_0$ in (\ref{M_gamma}) is known as the strong maximal operator defined on Heisenberg groups. The $\L^p$-boundedness of $\M_0$ has been proved by M. Christ \cite{Michael Christ 2}  using  a mixture of techniques developed previously by Ricci and Stein \cite{Ricci-Stein 1}-\cite{Ricci-Stein 2} and Christ \cite{Michael Christ 1} for singular integrals defined on sub-manifolds within a general setting of Nilpotent Lie groups.
\v
 {\bf Christ theorem }~~{\it Let $\M_0$ defined in  (\ref{M_gamma}). We have
\bel{Result Two}
\begin{array}{cc}\ds
\left\| \M_0 f\right\|_{\L^p(\R^{2n+1})}~\leq~\C_{p}~\left\| f\right\|_{\L^p(\R^{2n+1})},\qquad 1<p<\infty.
\end{array}
\eeq}\\
The remaining paper is organized as follows. In section 2, we show that two constraints consisting of   $\alpha,\beta$,$\vartheta$, $p,q$  in (\ref{Result One}) are necessary conditions. In section 3, we prove their sufficiency and finish the proof of  {\bf Theorem One}. In section 4, we prove {\bf Theorem Two}.

\section{Some necessary constraints}
\setcounter{equation}{0}
Let $\I_{\alpha\beta\vartheta}$ defined in (\ref{V})-(\ref{I}). By changing variable $\tau\mt \tau-\mu\big(u\cdot\eta-v\cdot\xi\big)$, we find
\bel{I rewrite}
\begin{array}{lr}\ds
\I_{\alpha\beta\vartheta} f(u,v,t)~=~\iiint_{\R^{2n+1}} f\left(\xi,\eta,\tau\right) \V^{\alpha\beta\vartheta}\Big[u-\xi, v-\eta, t-\tau-\mu(u\cdot\eta-v\cdot\xi)\Big]d\xi d\eta d\tau
\\\\ \ds~~~~~~~~~~~~~~~~~~~~
~=~\iiint_{\R^{2n+1}} f\left(\xi,\eta,\tau-\mu\big(u\cdot\eta-v\cdot\xi\big)\right) \V^{\alpha\beta\vartheta}(u-\xi,v-\eta,t-\tau) d\xi d\eta d\tau
\\\\ \ds~~~~~~~~~~~~~~~~~~~~
~=~\iiint_{\R^{2n+1}} f\left(\xi,\eta,\tau-\mu\big(u\cdot\eta-v\cdot\xi\big)\right) 
\\ \ds~~~~~~~~~~~~~~~~~~~~~~~~~~~~
|u-\xi|^{\alpha-n}|v-\eta|^{\alpha-n}|t-\tau|^{\beta-1} \Bigg[ {|u-\xi||v-\eta|\over |t-\tau|}+{|t-\tau|\over |u-\xi||v-\eta|}\Bigg]^{-\vartheta}d\xi d\eta d\tau.
\end{array}
\eeq
Consider a more general situation by replacing $\V^{\alpha\beta\vartheta}(u,v,t)$  with
\bel{V general}
|u|^{\alpha_1-n}|v|^{\alpha_2-n}|t|^{\beta-1} \Bigg[ {|u||v|\over |t|}+{|t|\over |u||v|}\Bigg]^{-\vartheta},\qquad \alpha_1,\alpha_2,\beta\in\R,\qquad \vartheta\ge0
\eeq
for $u\neq0$, $v\neq0$, $t\neq0$.

Let $(u,v,t)\mt (r u, s v, rs\lambda t)$ and $(\xi,\eta,\tau)\mt (r \xi, s \eta, rs\lambda \tau)$ for $r,s>0$ and  $0<\lambda<1$ or $\lambda>1$.  
We have
\bel{Dila I general}
\begin{array}{lr}\ds
\left\{ \iiint_{\R^{2n+1}} 
\left\{\iiint_{\R^{2n+1}} f\left[r^{-1} \xi,s^{-1}\eta,r^{-1}s^{-1}\lambda^{-1}\big[\tau-\mu\lambda\big(u\cdot\eta-v\cdot\xi\big)\big]\right] \right.\right.
\\ \ds
\left. \left.|u-\xi|^{\alpha_1-n}|v-\eta|^{\alpha_2-n}|t-\tau|^{\beta-1} \Bigg[ {|u-\xi||v-\eta|\over |t-\tau|}+{|t-\tau|\over |u-\xi||v-\eta|}\Bigg]^{-\vartheta}
 d\xi d\eta d\tau\right\}^q dudvdt\right\}^{1\over q}
\\\\ \ds
=~r^{\alpha_1+\beta}s^{\alpha_2+\beta}r^{n+1\over q}s^{n+1\over q}\lambda^\beta\lambda^{1\over q}~\left\{ \iiint_{\R^{2n+1}} 
\left\{\iiint_{\R^{2n+1}} f\left( \xi,\eta,\tau-\mu\big(u\cdot\eta-v\cdot\xi\big)\right) \right.\right.
\\ \ds~~~
\left. \left.|u-\xi|^{\alpha_1-n}|v-\eta|^{\alpha_2-n}|t-\tau|^{\beta-1} \Bigg[ {|u-\xi||v-\eta|\over \lambda|t-\tau|}+{\lambda|t-\tau|\over |u-\xi||v-\eta|}\Bigg]^{-\vartheta}
 d\xi d\eta d\tau\right\}^q dudvdt\right\}^{1\over q}
\\\\ \ds
\ge~r^{\alpha_1+\beta}s^{\alpha_2+\beta}r^{n+1\over q}s^{n+1\over q}\lambda^\beta\lambda^{1\over q}\left\{ \begin{array}{lr}\ds \lambda^\vartheta,\qquad0<\lambda<1,
\\ \ds
\lambda^{-\vartheta},\qquad \lambda>1
\end{array}\right.
\\ \ds~~~
\left\{ \iiint_{\R^{2n+1}} 
\left\{\iiint_{\R^{2n+1}} f\left( \xi,\eta,\tau-\mu\big(u\cdot\eta-v\cdot\xi\big)\right) \right.\right.
\\ \ds~~~
\left. \left.|u-\xi|^{\alpha_1-n}|v-\eta|^{\alpha_2-n}|t-\tau|^{\beta-1} \Bigg[ {|u-\xi||v-\eta|\over |t-\tau|}+{|t-\tau|\over |u-\xi||v-\eta|}\Bigg]^{-\vartheta}
 d\xi d\eta d\tau\right\}^q dudvdt\right\}^{1\over q}.
\end{array}
\eeq
The $\L^p\mt\L^q$-norm inequality in (\ref{Result Two}) implies that the last line of (\ref{Dila I general}) is bounded by
\bel{I L^p}
\begin{array}{lr}\ds
\left\{ \iiint_{\R^{2n+1}} \Big| f\left(r^{-1} \xi,s^{-1}\eta,r^{-1}s^{-1}\lambda^{-1} \tau\right) \Big|^p d\xi d\eta d\tau\right\}^{1\over p}~=~r^{n+1\over p}s^{n+1\over p} \lambda^{1\over p} \left\| f\right\|_{\L^p(\R^{2n+1})}.
\end{array}
\eeq
This must be true for every $r,s>0$ and $0<\lambda<1$ or $\lambda>1$. We necessarily have
\bel{Constraints alphabeta}
\begin{array}{cc}\ds
{\alpha_1+\beta\over n+1}~=~{1\over p}-{1\over q}~=~{\alpha_2+\beta\over n+1},
\qquad
\beta+\vartheta~\ge~{1\over p}-{1\over q}\qquad\hbox{\small{or}}\qquad \beta-\vartheta~\leq~{1\over p}-{1\over q}.
\end{array}
\eeq
The first constraint in (\ref{Constraints alphabeta}) forces us to have $\alpha_1=\alpha_2$. Therefore, write
\bel{Formula 12}
{\alpha+\beta\over n+1}~=~{1\over p}-{1\over q},\qquad \alpha=\alpha_1=\alpha_2.
\eeq
By bringing (\ref{Formula 12}) to the two inequalities in (\ref{Constraints alphabeta}), we find
\bel{theta computa}
\begin{array}{lr}\ds
\vartheta~\ge~\beta-{\alpha+\beta\over n+1}~=~{n\beta-\alpha\over n+1}\qquad \hbox{\small{or}}
\qquad
\vartheta~\ge~{\alpha+\beta\over n+1}-\beta~=~{\alpha-n\beta\over n+1}.
\end{array}
\eeq
Together, we conclude
$\vartheta~\ge~{|\alpha-n\beta|\over n+1}$.

\section{Proof of Theorem One}
\setcounter{equation}{0}
 Given $\alpha,\beta\in\R$ and $\vartheta\ge{|\alpha-n\beta |\over n+1}$, $\V^{\alpha\beta\vartheta}$ is a distribution in $\R^{2n+1}$ agree with $\V^{\alpha\beta\vartheta}(u,v,t)$ in (\ref{V}) whenever $u\neq0,v\neq0,t\neq0$.

Suppose  $\alpha\ge n\beta $. We have ${|\alpha-n\beta |\over n+1}={\alpha-n\beta \over n+1}$ and
\bel{EST1}
\begin{array}{lr}\ds
\V^{\alpha\beta\vartheta}(u,v,t)~\leq~|u|^{\alpha-n}|v|^{\alpha-n}|t|^{\beta-1} \Bigg[ {|u||v|\over |t|}+{|t|\over |u||v|}\Bigg]^{-{\alpha-n\beta \over n+1}}
\\\\ \ds~~~~~~~~~~~~~~~~~~~
~\leq~|u|^{\alpha-n}|v|^{\alpha-n}|t|^{\beta-1} \Bigg[ {|u||v|\over |t|}\Bigg]^{-{\alpha-n\beta \over n+1}}
\\\\ \ds~~~~~~~~~~~~~~~~~~~
~=~|u|^{n\big[{\alpha+\beta\over n+1}\big]-n}|v|^{n\big[{\alpha+\beta\over n+1}\big]-n} |t|^{{\alpha+\beta\over n+1}-1},\qquad\hbox{\small{$u\neq0$, $v\neq0$, $t\neq0$}}.
\end{array}
\eeq
Suppose $\alpha\leq n\beta $. We find  ${|\alpha-n\beta |\over n+1}={n\beta -\alpha\over n+1}$ and
\bel{EST2}
\begin{array}{lr}\ds
\V^{\alpha\beta\vartheta}(u,v,t)~\leq~|u|^{\alpha-n}|v|^{\alpha-n}|t|^{\beta-1} \Bigg[ {|u||v|\over |t|}+{|t|\over |u||v|}\Bigg]^{{\alpha-n\beta \over n+1}}
\\\\ \ds~~~~~~~~~~~~~~~~~~
~\leq~|u|^{\alpha-n}|v|^{\alpha-n}|t|^{\beta-1} \Bigg[ {|t|\over |u||v|}\Bigg]^{\alpha-n\beta \over n+1}
\\\\ \ds~~~~~~~~~~~~~~~~~~
~=~|u|^{n\big[{\alpha+\beta\over n+1}\big]-n}|v|^{n\big[{\alpha+\beta\over n+1}\big]-n} |t|^{{\alpha+\beta\over n+1}-1},\qquad \hbox{\small{$u\neq0$, $v\neq0$, $t\neq0$}}.
\end{array}
\eeq
Let $\I_{\alpha\beta\vartheta} f$ defined in (\ref{V})-(\ref{I}) and
\bel{Formula}
{\alpha+\beta\over n+1}~=~{1\over p}-{1\over q},\qquad 1<p<q<\infty.
\eeq
By changing variable $\tau\mt \tau-\mu\big(u\cdot\eta-v\cdot\xi\big)$, we have 
 \bel{I < =}
\begin{array}{lr}\ds
\I_{\alpha\beta\vartheta} f(u,v,t)
~=~\iiint_{\R^{2n+1}} f\left(\xi,\eta,\tau-\mu\big(u\cdot\eta-v\cdot\xi\big)\right) 
\\ \ds~~~~~~~~~~~~~~~~~~~~~~~~~~
|u-\xi|^{\alpha-n}|v-\eta|^{\alpha-n}|t-\tau|^{\beta-1} \Bigg[ {|u-\xi||v-\eta|\over |t-\tau|}+{|t-\tau|\over |u-\xi||v-\eta|}\Bigg]^{-\vartheta}d\xi d\eta d\tau
\\\\ \ds~~~~~~~~~~~~~~~~~~
~\leq~
\iiint_{\R^{2n+1}}  f\left(\xi,\eta,\tau-\mu\big(u\cdot\eta-v\cdot\xi\big)\right) 
\\ \ds~~~~~~~~~~~~~~~~~~~~~~~~~~
|u-\xi|^{n\big[{\alpha+\beta\over n+1}\big]-n}|v-\eta|^{n\big[{\alpha+\beta\over n+1}\big]-n} |t-\tau|^{{\alpha+\beta\over n+1}-1}d\xi d\eta d\tau\qquad\hbox{\small{by (\ref{EST1})-(\ref{EST2})}}.
\end{array}
\eeq
Because  $\V^{\alpha\beta\vartheta}$  is positive definite, it is suffice to assert $f\ge0$.
Define
\bel{F_beta}
\F_{\alpha\beta}(\xi,\eta, u,v,t)~=~\int_\R f\left(\xi,\eta,\tau-\mu\left(u \cdot\eta-v \cdot\xi\right)\right) |t-\tau|^{{\alpha+\beta\over n+1}-1} d\tau.
\eeq
From (\ref{I < =})-(\ref{F_beta}), we find
\bel{I rewrite F}
\begin{array}{lr}\ds
\I_{\alpha\beta\vartheta} f(u,v,t)
~\leq~\iint_{\R^{2n}} |u-\xi|^{n\big[{\alpha+\beta\over n+1}\big]-n}|v-\eta|^{n\big[{\alpha+\beta\over n+1}\big]-n} \F_{\alpha\beta}(\xi,\eta, u,v,t)d\xi d\eta.
\end{array}
\eeq
Recall  the {\bf Hardy-Littlewood-Sobolev theorem} stated in the beginning of this paper. By applying (\ref{HLS Ineq}) with $\a={\alpha+\beta\over n+1}$ and $\N=1$, we have
\bel{F regularity}
\begin{array}{lr}\ds
\left\{\int_{\R} \F^q_{\alpha\beta}(\xi,\eta, u,v,t) dt \right\}^{1\over q}~\leq~\C_{p~q} \left\{\int_{\R} \Big[ f\left(\xi,\eta,t-\mu\left(u \cdot\eta-v \cdot\xi\right)\right)\Big]^p dt\right\}^{1\over p}
\\\\ \ds~~~~~~~~~~~~~~~~~~~~~~~~~~~~~~~~~~~~~~~~~
~=~\C_{p~q} \left\| f(\xi,\eta,\cdot)\right\|_{\L^p(\R)}
\end{array}
\eeq
regardless of $(u,v)\in\R^n\times\R^n$.

On the other hand, by applying (\ref{HLS Ineq}) with $\a=n \big[{\alpha+\beta\over n+1}\big]$ and $\N=n$, we find
\bel{n Regularity u}
\begin{array}{lr}\ds
\left\{\int_{\R^n}         \left\{ \int_{\R^n}   |u-\xi|^{n \big[{\alpha+\beta\over n+1}\big]-n}  \left\| f(\xi,\eta,\cdot)\right\|_{\L^p(\R)} d\xi\right\}^q du\right\}^{1\over q}
~\leq~\C_{p~q} ~ \left\{\int_{\R^n} \left\| f(u,\eta,\cdot)\right\|_{\L^p(\R)}^p du\right\}^{1\over p}
\end{array}
\eeq
and
\bel{n Regularity v}
\begin{array}{lr}\ds
\left\{\int_{\R^n}         \left\{ \int_{\R^n}   |v-\eta|^{n \big[{\alpha+\beta\over n+1}\big]-n}  \left\| f(\xi,\eta,\cdot)\right\|_{\L^p(\R)} d\eta\right\}^q dv\right\}^{1\over q}
~\leq~\C_{p~q} ~ \left\{\int_{\R^n} \left\| f(\xi,v,\cdot)\right\|_{\L^p(\R)}^p dv\right\}^{1\over p}.
\end{array}
\eeq
From (\ref{I rewrite F}), we have
\bel{EST}
\begin{array}{lr}\ds
\left\|\I_{\alpha\beta\vartheta} f\right\|_{\L^q(\R^{2n+1})}
\\\\ \ds
~\leq~\left\{ \iiint_{\R^{2n+1}}\left\{\iint_{\R^{2n}} |u-\xi|^{n \big[{\alpha+\beta\over n+1}\big]-n}|v-\eta|^{n \big[{\alpha+\beta\over n+1}\big]-n} \F_{\alpha\beta}(\xi,\eta, u,v,t)d\xi d\eta\right\}^q dudvdt\right\}^{1\over q}
\\\\ \ds
~\leq~      \left\{\iint_{\R^{2n}}      \left\{ \iint_{\R^{2n}}   |u-\xi|^{n \big[{\alpha+\beta\over n+1}\big]-n}|v-\eta|^{n \big[{\alpha+\beta\over n+1}\big]-n}  \left\{\int_\R \F^q_{\alpha\beta}(\xi,\eta, u,v,t) dt \right\}^{1\over q} d\xi d\eta\right\}^q dudv\right\}^{1\over q}
\\ \ds~~~~~~~~~~~~~~~~~~~~~~~~~~~~~~~~~~~~~~~~~~~~~~~~~~~~~~~~~~~~~~~~~~~~~~~~~~~~~~~~~~~~~~~~~
\hbox{\small{by Minkowski integral inequality}}
\\\\ \ds
~\leq~\C_{p~q} \left\{\iint_{\R^{2n}}      \left\{ \iint_{\R^{2n}}   |u-\xi|^{n \big[{\alpha+\beta\over n+1}\big]-n}|v-\eta|^{n \big[{\alpha+\beta\over n+1}\big]-n}   \left\| f(\xi,\eta,\cdot)\right\|_{\L^p(\R)} d\xi d\eta\right\}^q dudv\right\}^{1\over q}\qquad\hbox{\small{by (\ref{F regularity})}}
\\\\ \ds
~\leq~\C_{p~q} \left\{\int_{\R^n}        \left\{\int_{\R^n}  \left\{ \int_{\R^n}   |v-\eta|^{n\big[{\alpha+\beta\over n+1}\big]-n}  \left\| f(u,\eta,\cdot)\right\|_{\L^p(\R)} d\eta\right\}^p du\right\}^{q\over p} dv\right\}^{1\over q}
\qquad
\hbox{\small{by  (\ref{n Regularity u})}}
\\\\ \ds
~\leq~\C_{p~q} \left\{\int_{\R^n} \left\{\int_{\R^n}         \left\{ \int_{\R^n}   |v-\eta|^{n\big[{\alpha+\beta\over n+1}\big]-n}  \left\| f(u,\eta,\cdot)\right\|_{\L^p(\R)} d\eta\right\}^q dv \right\}^{p\over q} du\right\}^{1\over p}
\\ \ds~~~~~~~~~~~~~~~~~~~~~~~~~~~~~~~~~~~~~~~~~~~~~~~~~~~~~~~~~~~
\hbox{\small{by Minkowski integral inequality}}
\\ \ds
~\leq~\C_{p~q}      \left\{\iint_{\R^{2n}}  \left\| f(u,v,\cdot)\right\|_{\L^p(\R)}^p du dv\right\}^{1\over p}\qquad\hbox{\small{by  (\ref{n Regularity v})}}
\\\\ \ds
~=~\C_{p~q} \left\| f\right\|_{\L^p(\R^{2n+1})}.
\end{array}
\eeq

\section{Proof of Theorem Two}
\setcounter{equation}{0}
Recall $\M_{\alpha\beta}$ defined in (\ref{M Z}). By taking $\xi\mt u-\xi$, $ \eta\mt v-\eta$ and $\tau\mt t-\tau$, it can be equivalently defined as
\bel{Mz}
\begin{array}{lr}\ds
\M_{\alpha\beta} f(u,v,t)~=~\sup_{\Rec\ni(u,v,t)\colon~\vol\{I\}=\vol\{\Q\}^{1\over n}\vol\{\P\}^{1\over n}}
\\ \ds~~~~~~~~~~~~~~~~~~~~~~~~~~~~
\vol\{\Q\}^{\frac{\alpha}{n}-1}\vol\{\P\}^{\frac{\alpha}{n}-1}\vol\{I\}^{\beta-1}\iiint_{\Rec}\left|f(\xi,\eta,\tau+\mu(u\cdot\eta-v\cdot\xi))\right|d\xi d\eta d\tau.
\end{array}
\eeq
Note that $\vol\{I\}=\vol\{\Q\}^{1\over n}\vol\{\P\}^{1\over n}$ implies
$\vol\{\Rec\}=\vol\{\Q\}\vol\{\P\}\vol\{I\}=\vol\{\Q\}^{1+{1\over n}} \vol\{\P\}^{1+ {1\over n}}$.
A direct computation shows
\bel{volume computa}
\begin{array}{lr}\ds
\vol\{\Q\}^{\frac{\alpha}{n}-1}\vol\{\P\}^{\frac{\alpha}{n}-1}\vol\{I\}^{\beta-1}~=~\vol\{\Q\}^{\big[\frac{\alpha+\beta}{n+1}-1\big]\left(1+\frac{1}{n}\right)}\vol\{\P\}^{\big[\frac{\alpha+\beta}{n+1}-1\big]\left(1+\frac{1}{n}\right)}
\\\\ \ds~~~~~~~~~~~~~~~~~~~~~~~~~~~~~~~~~~~~~~~~~~~~~~~~~
~=~\vol\{\Rec\}^{{\alpha+\beta\over n+1}-1}
\\\\ \ds~~~~~~~~~~~~~~~~~~~~~~~~~~~~~~~~~~~~~~~~~~~~~~~~~
~=~\vol\{\Q\}^{{\alpha+\beta\over n+1}-1} \vol\{\P\}^{{\alpha+\beta\over n+1}-1} \vol\{I\}^{{\alpha+\beta\over n+1}-1}.
\end{array}
\eeq
Observe that for every $\Rec=\Q\times\P\times I\ni(u,v,t)$ with $\vol\{I\}=\vol\{\Q\}^{1\over n}\vol\{\P\}^{1\over n}$, we have
\bel{Compare}
\begin{array}{lr}\ds
\vol\{\Q\}^{\frac{\alpha}{n}-1}\vol\{\P\}^{\frac{\alpha}{n}-1}\vol\{I\}^{\beta-1}\iiint_{\Rec\ni(u,v,t)}\left|f(\xi,\eta,\tau+\mu(u\cdot\eta-v\cdot\xi))\right|d\xi d\eta d\tau
\\\\ \ds
=~\vol\{\Q\}^{{\alpha+\beta\over n+1}-1} \vol\{\P\}^{{\alpha+\beta\over n+1}-1} \vol\{I\}^{{\alpha+\beta\over n+1}-1} 
\\ \ds~~~~
\iiint_{\Rec\ni(u,v,t)}\left|f(\xi,\eta,\tau+\mu(u\cdot\eta-v\cdot\xi))\right|d\xi d\eta d\tau\qquad\hbox{\small{by (\ref{volume computa})}}
\\\\ \ds
\leq~\iiint_{\R^{2n+1}} \left|f\left(\xi,\eta,\tau+\mu\big(u\cdot\eta-v\cdot\xi\big)\right) \right|
|u-\xi|^{n\big[{\alpha+\beta\over n+1}\big]-n}|v-\eta|^{n\big[{\alpha+\beta\over n+1}\big]-n}|t-\tau|^{{\alpha+\beta\over n+1}-1} d\xi d\eta d\tau.
\end{array}
\eeq
From (\ref{Mz}) and (\ref{Compare}), we conclude
\bel{M<I}
\begin{array}{lr}\ds
\M_{\alpha\beta} f(u,v,t)~\leq~
\\\\ \ds
\iiint_{\R^{2n+1}} \left|f\left(\xi,\eta,\tau+\mu\big(u\cdot\eta-v\cdot\xi\big)\right) \right|
|u-\xi|^{n\big[{\alpha+\beta\over n+1}\big]-n}|v-\eta|^{n\big[{\alpha+\beta\over n+1}\big]-n}|t-\tau|^{{\alpha+\beta\over n+1}-1} d\xi d\eta d\tau
\end{array}
\eeq
for every $(u,v,t)\in\R^n\times\R^n\times\R$.

Let 
\bel{formula}
{\alpha+\beta\over n+1}~=~{1\over p}-{1\over q},\qquad 1<p<q<\infty.
\eeq
By carrying out the same estimates in (\ref{F_beta})-(\ref{EST}), we obtain
\bel{Result Two again}
\left\| \M_{\alpha\beta} f\right\|_{\L^q(\R^{2n+1})}~\leq~\C_{p~q}~\left\| f\right\|_{\L^p(\R^{2n+1})}.
\eeq
Next, we show (\ref{formula})-(\ref{Result Two again}) also hold for $\alpha+\beta=0$ and $1<p=q<\infty$.

Recall $\M_\gamma$ defined in (\ref{M_gamma}) for $0\leq\gamma<1$. By taking $\xi\mt u-\xi$, $ \eta\mt v-\eta$ and $\tau\mt t-\tau$, $\M_\gamma$ can be equivalently defined as
\bel{M_gamma equi}
\M_\gamma f(u,v,t)~=~
\sup_{\Rec\ni(u,v,t)}\vol\{\Rec\}^{\gamma-1}\iiint_{\Rec}\left|f(\xi,\eta,\tau+\mu(u\cdot\eta-v\cdot\xi))\right|d\xi d\eta d\tau
\eeq
where $\Rec$ in (\ref{M_gamma equi}) is any rectangle ( possibly $2n+1$-parameter ) parallel to the coordinates.

Let $\gamma={\alpha+\beta\over n+1}$. From (\ref{Mz})-(\ref{volume computa}), we find
\bel{Mz<}
\begin{array}{lr}\ds
\M_{\alpha\beta} f(u,v,t)~=~\sup_{\Rec\ni(u,v,t)\colon~\vol\{I\}=\vol\{\Q\}^{1\over n}\vol\{\P\}^{1\over n}}
\\ \ds~~~~~~~~~~~~~~~~~~~~~~~~~~~~
\vol\{\Q\}^{\frac{\alpha}{n}-1}\vol\{\P\}^{\frac{\alpha}{n}-1}\vol\{I\}^{\beta-1}\iiint_{\Rec}\left|f(\xi,\eta,\tau+\mu(u\cdot\eta-v\cdot\xi))\right|d\xi d\eta d\tau
\\\\ \ds
~=~\sup_{\Rec\ni(u,v,t)\colon\vol\{I\}=\vol\{\Q\}^{1\over n}\vol\{\P\}^{1\over n}}\vol\{\Rec\}^{\frac{\alpha+\beta}{n+1}-1}\iiint_{\Rec}\left|f(\xi,\eta,\tau+\mu(u\cdot\eta-v\cdot\xi))\right|d\xi d\eta d\tau
\\\\ \ds
~\leq~\sup_{\Rec\ni(u,v,t)}\vol\{\Rec\}^{\gamma-1}\iiint_{\Rec}\left|f(\xi,\eta,\tau+\mu(u\cdot\eta-v\cdot\xi))\right|d\xi d\eta d\tau
\\\\ \ds
~=~\M_\gamma f(u,v,t).
\end{array}
\eeq
Hence, $\M_{\alpha\beta}$ is controlled by the strong fractional maximal operator $\M_\gamma$ whenever  $\gamma={\alpha+\beta\over n+1}$. At $\gamma=0$, {\bf Christ theorem} stated in Section 1 implies
\bel{M_alpha beta 0}
\begin{array}{cc}\ds
\left\| \M_{\alpha\beta} f\right\|_{\L^p(\R^{2n+1})}~\leq~\C_{p}~\left\| f\right\|_{\L^p(\R^{2n+1})},\qquad 1<p<\infty
\\\\ \ds
\hbox{whenever} \qquad \alpha+\beta~=~0.
\end{array}
\eeq

{\small Westlake University, Hangzhou, 310030, China}\\
{\small email: wangzipeng@westlake.edu.cn}

{\small School of Mathematical Sciences, Zhejiang University, Hangzhou, 310058, China}\\
{\small email: sunchuhan@zju.edu.cn}

\end{document}